\magnification = 1200
\overfullrule = 0pt
\parindent = 0pt

\def\title#1{\centerline{\bf#1}}
\def\frac#1#2{{#1\over #2}}%

\title{Another Proof of a Famous Inequality}
\bigskip
{\bf 1.}\quad The mean-value theorem of differential  calculus, in its simplest
form, states that if a function $\phi$ has a derivative at all points  between $c$ and
$d$, and including those points,  then for some $e$ {\it strictly} between
$c$ and $d$  
$$
\phi(d) - \phi(c) = (d-c)\phi'(e).
$$
 The word ``between''  was used as the order of $c$ and $d$ is unimportant. However in the 
following we will assume that $c<d$.
\smallskip
{\bf 2.}\quad The basic inequality of the title is: if $a\ne b$ are positive  and if $u,v>0$ then
 $$
 \frac{ua + vb}{u+v} > \bigl(a^ub^v\bigr)^{1/(u+v)}.\eqno(GA)
$$
This is the simplest case of the inequality between the arithmetic and geometric means, from
which the general case follows by a simple induction; see [2, p.81].

Two things should be remarked for our later discussion: 
\smallskip
\noindent (i) as both sides of (GA) are symmetric
in $a$ and $b$ there is no loss in generality in assuming $a<b$; 

\noindent (ii) as (G) is unchanged if $ u,v$ are replaced by $\lambda u, \lambda v$
respectively, for any $\lambda>0$, we can assume,  without loss in generality,  that $u,
 v \ge 1$
\smallskip
 {\bf 3.}\quad Consider then
the  distinct power 
functions
$\phi_1(x) = x^u,\;
\phi_2(x) = x^v$; where $u\ne v$, $u,v\ge 1$ and $x>0$ . Applying the mean value
theorem  to both of these functions we get, on cancelling the common factor $d-c$,
$$
 \frac{d^u- c^u}{d^v-c^v} = \frac{ue_1^{u-1}}{ve_2^{v-1}}\, ,\eqno(1)
$$
 for some $e_1, e_2$,  with $c<e_1, e_2<d$.

Using the various conditions   (1) leads to
$$
 \frac{d^u- c^u}{d^v-c^v} > \frac{uc^{u-1}}{vd^{v-1}}>\frac{uc^u}{vd^v};\eqno(2)
$$
 or on  multiplying out

$$
 vd^v(d^u-c^u)>uc^u(d^v-c^v)\eqno(3)
$$
 To get (3) we have assumed that $u\ne v$ but  it is easy to check that (3)
remains valid when $u=v$.

Rewriting (3) gives
$$
uc^{u+v} + vd^{u+v} > (u+v) c^ud^v.\eqno(4)
$$
Inequality (4) is the just (GA) as  can be
seen by putting $a= c^{u+v},\; b= d^{u+v}$.
\smallskip
{\bf 3.}\quad  It is possible
to have $e_1 = e_2 $ in (1). However then we must use the more
sophisticated Cauchy mean-value theorem: if $\phi_1,\phi_2$ satisfy the
condition that $\phi$ satisfies section 1 above, and if in addition $ \phi_2'$ is
never zero then for some $e$ strictly between $c$ and $d$,
 $$
\frac{\phi_1(d) -\phi_1(c)}{\phi_2(d) - \phi_2(c)} =
\frac{\phi_1'(e)}{\phi_2'(e)}\, . 
$$
 Now choosing the $\phi_1,\phi_2$ as in section 2 above we get (1) with $ e_1
= e_2$, that is 
$$
 \frac{d^u- c^u}{d^v-c^v} = \frac{ue^{u-1}}{ve^{v-1}} =\frac{ue^u}{ve^v}.
$$
The only advantage of this approach is that we now need not assume $ u,v>1$, $ u,v>0$ will
suffice.
\smallskip
{\bf 4.}\quad This proof of  (GA) is based on an idea in [1]; see also [2,p. 78]. There $u, v$ 
are integers bigger than $1$ and to get (2) the identity 
$$
\frac{ d^n - c^n }{d-c}= c^{n-1}+ c^{n-2}d +\cdots + cd^{n-2} + d^{n-1},
$$
obtained by summing a geometric series, is used. For then in the case $c<d$ and
$u,v\ge 2$, integers, 
$$
\frac{d^u- c^u}{d^v-c^v}=\frac{ c^{u-1} +\cdots + d^{u-1}}{
c^{v-1}+\cdots + d^{v-1}} >\frac{uc^{u-1}}{vd^{v-1}}\, ,
$$
 which is (2).
\smallskip
\centerline{\bf Bibliography}
\smallskip
1. T.Arnold Brown. Elementary inequalities, {\sl\ Math.\ Gazette}, \enspace 263\
(1941),\ 2--11.

2. P.S.Bullen.  {\sl Handbook of Means and Their Inequalities\/},   Kluwer Academic  Publishers,
Dordrecht/Boston/London.\ 2003.

\bye